\documentclass{amsart}

\usepackage{amsmath,amssymb,amsbsy,amsthm,amsfonts,mathtools,dsfont,hyperref,mathrsfs}

\usepackage{enumitem}
\usepackage{wrapfig}
\usepackage[top=3.5cm,bottom=3.5cm,left=3.5cm,right=3.5cm]{geometry}
\usepackage[normalem]{ulem}

\usepackage{graphicx}
\usepackage{tikz}
\usepackage{tikzscale}
\usetikzlibrary{intersections,decorations.markings}
\tikzset{>=latex}
\usepackage{pgfplots}
\pgfplotsset{compat=newest}
\usepgfplotslibrary{fillbetween}
\usepackage[normalem]{ulem}

\allowdisplaybreaks[3]

\newcommand{\tr}{\operatorname{Tr}}

\begin{document}

\newcommand{\ddd}{\,{\rm d}}

\def\note#1{\marginpar{\small #1}}
\def\tens#1{\pmb{\mathsf{#1}}}
\def\vec#1{\boldsymbol{#1}}
\def\norm#1{\left|\!\left| #1 \right|\!\right|}
\def\fnorm#1{|\!| #1 |\!|}
\def\abs#1{\left| #1 \right|}
\def\ti{\text{I}}
\def\tii{\text{I\!I}}
\def\tiii{\text{I\!I\!I}}

\newcommand{\loc}{{\rm loc}}
\def\diver{\mathop{\mathrm{div}}\nolimits}
\def\grad{\mathop{\mathrm{grad}}\nolimits}
\def\Div{\mathop{\mathrm{Div}}\nolimits}
\def\Grad{\mathop{\mathrm{Grad}}\nolimits}
\def\cof{\mathop{\mathrm{cof}}\nolimits}
\def\det{\mathop{\mathrm{det}}\nolimits}
\def\lin{\mathop{\mathrm{span}}\nolimits}
\def\pr{\noindent \textbf{Proof: }}

\def\pp#1#2{\frac{\partial #1}{\partial #2}}
\def\dd#1#2{\frac{\d #1}{\d #2}}
\def\vec#1{\boldsymbol{#1}}

\def\0{\vec{0}}
\def\A{\mathcal{A}}
\def\B{\mathcal{B}}
\def\b{\vec{b}}
\def\C{\mathcal{C}}
\def\c{\vec{c}}
\def\D{\vec{Dv}}
\def\DD{\vec{D}}
\def\BB{\vec{B}}
\def\e{\varepsilon}
\def\er{\epsilon}
\def\f{\vec{f}}
\def\F{\vec{F}}
\def\tF{\tilde{\F}}
\def\g{\vec{g}}
\def\G{\vec{G}}
\def\cG{\mathcal{\G}}
\def\H{\vec{H}}
\def\cH{\mathcal{H}}
\def\I{\vec{I}}
\def\Im{\text{Im}}
\def\j{\vec{j}}
\def\J{\vec{J}}
\def\dd{\vec{d}}
\def\k{\vec{k}}
\def\n{\vec{n}}
\def\q{\vec{q}}
\def\S{\vec{S}}
\def\s{\vec{s}}
\def\T{\vec{T}}
\def\u{\vec{u}}
\def\vp{\vec{\varphi}}
\def\vv{\vec{v}}
\def\vvt{\vv_{\tau}}
\def\vov{\vv\otimes\vv}
\def\cV{\mathcal{V}}
\def\w{\vec{w}}
\def\W{\vec{W}}
\def\x{\vec{x}}
\def\z{\vec{Z}}
\def\tz{\tilde{\z}}
\def\Z{\vec{Z}}
\def\X{\vec{X}}
\def\Y{\vec{Y}}
\def\balfa{\vec{\alpha}}

\def\Ge{\G_{\e}}
\def\ge{\g_{\e}}
\def\fidv{\phi_{\delta}(|\vv|^2)}
\def\fidve{\phi_{\delta}(|\ve|^2)}
\def\fidvd{\phi_{\delta}(|\vd|^2)}

\def\Ae{\A_\e}
\def\Aee{\Ae^\e}
\def\Aeetilde{\tilde{A}_\e^\e}
\def\Be{\B_\e}
\def\Bee{\Be^\e}
\def\De{\DD\ve}
\def\DDe{\DD^\e}
\def\Dvd{\DD\vd}
\def\oD{\overline{\DD}}
\def\tD{\tilde{\DD}}
\def\Dn{\DD^\e}
\def\Dno{\overline{\Dn}}
\def\Dnt{\tilde{\Dn}}
\def\Dm{\DD^\eta}
\def\Dmo{\overline{\Dm}}
\def\Dmt{\tilde{\Dm}}
\def\Se{\S^\e}
\def\se{\s^\e}
\def\ose{\overline{\se}}
\def\oS{\overline{\S}}
\def\tS{\tilde{\S}}
\def\Sn{\S^\e}
\def\Sno{\overline{\Sn}}
\def\Snt{\tilde{\Sn}}
\def\Sm{\S^\eta}
\def\Smo{\overline{\Sm}}
\def\Smt{\tilde{\Sm}}
\def\ve{\vv^\e}
\def\ove{\overline{\ve}}
\def\vove{\ve\otimes\ve}
\def\vd{\vv^\delta}
\def\sd{\s^\delta}
\def\Sd{\S^\delta}
\def\Dd{\DD\vd}

\def\Wnd#1{W^{1,#1}_{\n, \diver}}
\def\Wndr{W^{1,r}_{\n, \diver}}

\def\o{\Omega}
\def\po{\partial \Omega}
\def\dt{\frac{\d}{\d t}}
\def\pt{\partial_t}
\def\it{\int_0^t \!}
\def\iT{\int_0^T \!}
\def\io{\int_{\o} \!}
\def\iq{\int_{Q} \!}
\def\iqt{\int_{Q^t} \!}
\def\ipo{\int_{\po} \!}
\def\ig{\int_{\Gamma} \!}
\def\igt{\int_{\Gamma^t} \!}

\def\d{\, \textrm{d}}

\def\mn{\mathcal{P}}
\def\du{\mathcal{W}}
\def\tr{\operatorname{tr}}
\def\tow{\rightharpoonup}


\newtheorem{theorem}{Theorem}[section]
\newtheorem{lemma}[theorem]{Lemma}
\newtheorem{proposition}[theorem]{Proposition}
\newtheorem{remark}[theorem]{Remark}
\newtheorem{corollary}[theorem]{Corollary}
\newtheorem{definition}[theorem]{Definition}
\newtheorem{example}[theorem]{Example}

\numberwithin{equation}{section}

\title[Generalized Korteweg's fluids]{A generalization of the classical Euler and Korteweg fluids}

\thanks{K.~R.~Rajagopal thanks the Office of Naval Research for its support of this work.}

\author[K.~R.~Rajagopal]{K.~R.~Rajagopal}
\address{Department of Mechanical Engineering,  
Texas A\&M University, College Station, TX 77845 USA}
\email{krajagopal@tamu.edu}


\begin{abstract} 
The aim of this short paper is threefold. First, we develop an implicit generalization of a constitutive relation introduced by Korteweg \cite{korteweg.dj:sur} that can describe the phenomenon of capillarity. Second, using a sub-class of the constitutive relations (implicit Euler equations), we show that even in that simple situation more than one of the members of the sub-class may be able to describe one or a set of experiments one is interested in describing, and we must determine which amongst these constitutive relations is the best by culling the class by systematically comparing against an increasing set of observations. (The implicit generalization developed in this paper is not a sub-class of the implicit generalization of the Navier-Stokes fluid developed by Rajagopal \cite{krr1, krr2} or the generalization due to Pr\accent23u\v{s}a and Rajagopal \cite{prusa.rajagopal.kr:on}, as spatial gradients of the density appear in the constitutive relation developed by Korteweg \cite{korteweg.dj:sur}.) Third, we  introduce a challenging set of partial differential equations that would lead to new techniques in both analysis and numerical analysis to study such equations.
\end{abstract}

\maketitle

\emph{On the occasion of Josef M\' alek’s sixtieth birthday, steadfast friend, trusty collaborator.} 

\bigskip

\section{Introduction} \label{Sec1}

Rajagopal \cite{krr1} introduced algebraic implicit constitutive relations to describe the response of both solids and fluids, and later Pr\accent23u\v{s}a and Rajagopal \cite{prusa.rajagopal.kr:on} generalized the class of simple materials introduced by Noll \cite{noll.w:mathematical} to the class of implicit constitutive relations between the history of the stress and the history of the deformation gradient, and showed that under the assumption of fading memory, when the appropriate approximations are carried out, the implicit relations yield  both differential type and rate type approximations. The approximations obtained by Coleman and Noll  \cite{coleman.bd.noll.w:approximation} within the context of Simple fluids are a special sub-set of the approximations obtained by Pr\accent23u\v{s}a and Rajagopal \cite{prusa.rajagopal.kr:on}.
This study was followed by the generalization by Rajagopal \cite{krr15} who studied the anisotropy of implicit constitutive relation between the histories of the density, stress and the deformation gradient. All the above implicit constitutive relations do not include spatial gradient of the density, and thus the constitutive relation introduced by Korteweg \cite{korteweg.dj:sur} is not a sub-class of these implicit constitutive relations introduce by Rajagopal \cite{krr1}, Pr\accent23u\v{s}a and Rajagopal \cite{prusa.rajagopal.kr:on} and Rajagopal \cite{krr15}.

Thus, the first objective of the short work is to develop implicit generalizations of the classical Korteweg fluids (which includes implicit generalizations of the Euler fluid), with a view towards increasing the arsenal of the modeler to describe the response of compressible fluids. The second objective addresses the issue of determining constitutive relation that can describe observed phenomena, a problem confronted by the modeler. Conjectures are propounded on the basis of observations, and the iterative process between carefully constructed experiments to test these conjectures. The back and forth between the polishing of conjectures and refining of experiments hopefully leads to a theory that is simple, economical, with predictive capability, allowing for consilience of inductions and falsifiability. Most “theories” that are in vogue do not rise to such levels; they merely explain a small body of evidence. A constitutive theory is an explanation for the response of a particular class of materials\footnote{The terminology “constitutive theory” or “constitutive relations” is a misnomer. Constitutive applies to how a material is constituted (how it is composed), but it is used as a synonym for response functions (see Rajagopal \cite{krr23} for a detailed discussion of this erroneous usage of the terminology “constitutive relations”).} , based on our observation of how these materials behave when subject to external stimuli. The question then arises, given a class of flows that have been observed of a particular fluid, what are the class of constitutive relation that best explain the class of observed flows\footnote{An interesting mathematical generalization of this question, within the context of ordinary differential or partial differential equations is the following: given a particular class of solutions determine the class of ordinary differential or partial differential equations wherein such a class of solutions is possible? Within the context of incompressible isotropic Green elasticity (see Green \cite{green1837}, Truesdell and Noll \cite{ truesdell.c.noll.w:non-linear*1}) such a question has been investigated by McLeod, Rajagopal and Wineman \cite{McLeodkrr1988} for a class of inhomogeneous shear superposed on homogeneous triaxial extension. They delineate a class of stored energy functions which lead to a class of ordinary differential equations for which they prove solutions exist.}. Using the class of implicit Euler fluids, and a very simple static solution, we show that infinity of constitutive relations could describe the solution. This large class has to then be whittled down by considering more and more observed flows, arriving at a reasonable constitutive relation. Finally, we remark on the system of partial differential equations that arise from the implicit constitutive relations that we develop which can be the food for thought to mathematical and numerical analysts.

In a series of papers published between 1754 and 1761, Euler \cite{euler1, euler2, euler3} 
developed an idealized fluid model that has proved extremely useful in describing the flows of a large class of fluids. The Euler fluid\footnote{In the first paper on inviscid fluids, Euler required that the vorticity be zero in flows of incompressible inviscid fluids. Later, he generalized the investigation to include both compressible and incompressible fluid and relaxed the requirement that the flows be irrotational.}  is defined by a constitutive expression for the stress in terms of the density, namely (see Truesdell \cite{truesdell1977}, Truesdell and Rajagopal \cite{truesdell.c.rajagopal.kr:introduction})
\begin{equation}
    \label{1}
    \T = -p(\rho)\I,
\end{equation}
where $\rho$ is the density and $p$ is the mean value of the stress referred to as the mechanical pressure, and $\I$ denotes the identity tensor\footnote{The incompressible counterpart of the constitutive relation \eqref{1} takes the form $\T = -p\I$, where $p$ is the indeterminate scalar that is a consequence of the constraint of incompressibility.}. The Euler fluid is a perfectly elastic fluid incapable of dissipation. We shall see later that \eqref{1} is a very special sub-class of the Korteweg fluid whose constitutive relation is given by \eqref{4}. Usually, one also allows for the effect of temperature in which case the constitutive expression takes the form
\begin{equation}
    \label{2}
    \T = -p(\rho,\theta)\I,
\end{equation}
where now $\theta$ is the temperature. The classical ideal gas is an example of a Euler fluid. In \eqref{2}, $p$ is referred to as the thermodynamic pressure. 

The expression for the thermodynamic pressure $p$ as a function of the density $\rho$ and the temperature $\theta$, which is usually referred to as the equation of state, relates the various quantities that appear in it. In classical thermodynamics, one seems to take the approach that the quantities that appear in the equation of state are related and there does not seem much deliberation with regard to which of these quantities might be a cause and which an effect. One finds the thermodynamic pressure expressed in terms of the density and temperature; the density being expressed in terms of the thermodynamic pressure and temperature; or the temperature being expressed in terms of the thermodynamic pressure and density. While the notion of temperature is not a primitive when one is dealing within the context of statistical thermodynamics, it is presumed to be the cause for the transfer of “heat” (energy in thermal form) within the context of classical thermodynamics. For instance, Maxwell \cite{maxwell1871} states through \emph{``The temperature of a body is its thermal state considered with reference to power of communicating heat to other bodies."} That is, temperature is the power to transmit heat. Fosdick and Rajagopal \cite{Fosdickkrr83} have shown that the notion of transfer of heat (transfer of thermal energy) implies the existence of a locally Euclidean Hausdorff space of one dimension, namely the existence of temperature, that is the existence of the concept of temperature is a necessary precursor for heat transfer to take place, and it is the difference in temperature that causes heat transfer to take place, which is usually described by Fourier’s law: $\q=-k \nabla \theta$, where $\q$ is the heat flux (the effect), $\theta$ is the temperature (the gradient of $\theta$  is the cause for the heat flux) and $k$ is the thermal conductivity.

In the case of the constitutive representation \eqref{2}, the mechanical pressure (mean value of the stress) and thermodynamic pressure are the same. This is not always the case, especially when one considers the compressible Navier-Stokes fluid, many non-Newtonian fluids or the Korteweg fluid. The indiscriminate use of the terminology “pressure” has been the cause for much confusion (see Rajagopal \cite{krr5}) as will become clear from what follows. 

Korteweg \cite{korteweg.dj:sur} developed a constitutive relation for a fluid wherein the stress depended on both the density, its first and second spatial gradients, and the symmetric part of the velocity gradient. Models wherein the stress depends on the density and the gradients of the density have also been used to describe the response of granular materials (see Goodman and Cowin \cite{Goodman1972}, Hutter and Rajagopal \cite{hutterkrr1994}). The important fact to bear in mind is that the gradient in question is the Eulerian spatial gradient and the constitutive models under consideration are models for homogeneous bodies. M\' alek and Rajagopal \cite{Makrr06} have looked at constitutive equations for inhomogeneous bodies wherein they considered the possibility of the stress depending on the Lagrangian spatial gradient of the density. Here, we shall only consider the implications of constitutive relations for homogeneous bodies. While we shall not consider inhomogeneous bodies in this short paper, it is easy to generalize the results established here to the case of inhomogeneous bodies.

Korteweg \cite{korteweg.dj:sur} proposed a constitutive expression which falls into the class of materials defined through
\begin{equation}
    \label{3}
    \T = \vec{f}(\rho,\theta, \nabla \rho, \nabla^{(2)}\rho, \D).
\end{equation}
The Korteweg fluid takes the special form\footnote{Notice that when $\D = 0$, that is when there is no flow, and when we set $p(\rho)=\alpha_0(\rho) + \alpha_1 \Delta \rho + \alpha_2 |\nabla\rho|^2$, the mean value of the stress for the constitutive relation \eqref{4} is not necessarily the pressure $p(\rho)$, as there being contributions from the other terms.}:
\begin{equation}
    \label{4}
    \T = \left(\alpha_0(\rho) + \alpha_1 \Delta \rho + \alpha_2 |\nabla\rho|^2\right)\I + \alpha_3 (\nabla\rho \otimes\nabla \rho) + \alpha_4 \nabla ^{(2)}\rho + \lambda \tr \D + 2\mu \D,
\end{equation}
where $\alpha_0$ is a function depending on the density and $\alpha_i$, $i=1,2,3,4$, $\mu$ and $\lambda$ are constants and $\D$ is the symmetric part of the velocity gradient.  The above model presents lots of challenges with regard to the solution of boundary value problems as the balance of linear momentum will in general contain three spatial derivatives of the density and it is far from apparent the three boundary conditions that ought to be prescribed. Recently, Sou\v{c}ek, Heida and M\'{a}lek \cite{soucek2020} have determined boundary conditions for the Korteweg like fluids on the basis of thermodynamics. The boundary condition for the density is given in terms of the normal derivatives of the density at the boundary, but this condition might not be useful in the resolution of some general boundary value problems. We notice that to describe the fluid using the constitutive relation \eqref{4}, we have to be able to prescribe seven material constants, and it is a tremendously daunting task to delineate an experimental program in which these constants can be measured\footnote{Even within the context of the classical Navier-Stokes fluid, while one can measure the shear viscosity, one cannot determine the other viscosity that appears in the representation for the stress. One can measure the bulk modulus $\kappa=3\lambda+2\mu$, and thus indirectly compute $\lambda$.}.

The models considered by Euler and Korteweg are explicit expressions for the stress in terms of the dependent variables. Recently, Rajagopal \cite{krr1,krr2,krr3} introduced implicit constitutive relations that relate the stress and an appropriate kinematical quantity depending on whether the model describes the response of a solid or a fluid. In the case of fluid models, Rajagopal \cite{krr1,krr2} introduced constitutive relations where the stress, the density and the symmetric part of the velocity gradient were related through an algebraic relation. M\'{a}lek et al. \cite{mpr10} studied generalizations of the classical Navier-Stokes constitutive equation within the class of the implicit model introduced by Rajagopal, and Le Roux and Rajagopal developed implicit relations which in a simple shear flow exhibited a non-monotone response between the shear rate and shear stress. Perl\'{a}cov\'{a} and Pr\accent23u\v{s}a \cite{PerPr} showed that sub-classes of the implicit constitutive relations developed by Le Roux and Rajagopal \cite{LerouxKRR} can describe the response of colloidal solutions.  

In addition to implicit generalizations of the Euler and Korteweg fluids, it is possible to systematically develop implicit generalizations of many other constitutive relations including fluids with thresholds. Blechta, M\' alek and Rajagopal \cite{blechta2020} have provided a taxonomy of constitutive relations for incompressible fluids. It ought to be possible to develop such a categorization for compressible fluids as well (see also M\' alek and Rajagopal \cite{mr10}).

\section{Generalization of the classical Korteweg fluid model}\label{Sec2}

We notice that the second gradient of density appears in the Korteweg constitutive relation \eqref{4}. An implicit generalization of this is the implicit constitutive relation \eqref{3}. In what follows, we shall ignore the dependence of the constitutive relations on the second gradient of the density. Ignoring the dependence of the second gradient of the density implies that we cannot recover the full Korteweg model within the context of the implicit equation that is being proposed. However, inclusion of the second gradient of density and using standard representation theorems would lead to a constitutive relation characterized by a large a set of material functions that it would be impossible to fashion a meaningful experimental program to determine all of them, even if these material functions are assumed to be constants. Even the simplified model that we consider requires the determination of six material functions of the density, numerous invariants of tensors involving the stress, the tensor product of the gradient of the density with itself, and the product of the stress and second gradient of the density, as demanded by representation theorems (see Spencer \cite{spencer.ajm:theory}). Using such implicit models would require one to make sensible simplifications of the constitutive relations for them to be of any use whatsoever.

Let us consider an implicit relation between the Cauchy stress $\T$ the density $\rho$ and $\nabla \rho$, given by
\begin{equation}
    \label{5}
    \vec{f}(\rho,\nabla \rho, \T) = \vec{0}.
\end{equation}
If the fluid under consideration is isotropic, then $\vec{f}$ has to meet
\begin{equation}
    \label{6}
    \vec{f}(\rho,\vec{Q}\nabla \rho, \vec{Q}\T\vec{Q}^T) = \vec{Q} \vec{f}(\rho, \nabla\rho, \T) \vec{Q}^T  \qquad \textrm{ for all } \vec{Q}\in \mathcal{O},
\end{equation}
where $\mathcal{O}$ is the orthogonal group. Standard representation theorems (see Spencer \cite{spencer.ajm:theory}) then lead to the representation
\begin{equation}
    \label{7}\begin{split}
    \alpha_1 \I + \alpha_2 \T &+ \alpha_3 (\nabla\rho \otimes\nabla \rho) + \alpha_4 \T^2 \\ &+ \alpha_5 \left( \nabla \rho \otimes \T \nabla\rho + \T \nabla \rho \otimes\nabla\rho\right) + \alpha_6 \left( \nabla \rho \otimes \T^2 \nabla\rho + \T^2 \nabla \rho \otimes\nabla\rho\right) = \vec{0},\end{split}
\end{equation}
where the material moduli $\alpha_i$, $i=1, \dots, 6$, depend on the density $\rho$ and the following invariants:
\begin{equation}
    \label{8}
    \tr \T, \tr \T^2, \tr \T^3, \tr (\nabla \rho\otimes\nabla \rho), \tr (\nabla \rho\otimes\T \nabla \rho), \tr (\nabla\rho\otimes\T^2 \nabla \rho).
\end{equation}
If we restrict the implicit constitutive relation \eqref{7} to being linear in the Cauchy stress, we obtain
\begin{equation}
    \label{9}
    \alpha_1 \I + \alpha_2 \T + \alpha_3 (\nabla\rho \otimes\nabla \rho) + \alpha_5 \left( \nabla \rho \otimes \T \nabla\rho + \T \nabla \rho \otimes\nabla\rho\right) = \vec{0},
\end{equation}
where the material moduli $\alpha_1$ and $\alpha_3$ depend on the density $\rho$ and the following invariants 
\begin{equation}
    \label{10}
    \tr \T, \tr (\nabla\rho\otimes\nabla\rho), \tr(\nabla\rho\otimes\T \nabla\rho),
\end{equation}
and the material moduli $\alpha_2$ and $\alpha_5$ depend on the density $\rho$ and $\tr(\nabla\rho\otimes\nabla \rho)$. 

We now proceed to show that given an experimental observation, even within a much simpler sub-class of constitutive relations, implicit generalization of the Euler equation, infinite number of constitutive relations could explain a particular phenomenon, and that we have to whittle them down by considering several experimental observations.  Let us consider the sub-class of \eqref{5} wherein we have an implicit relation between the density and the stress. In this case
\begin{equation}
\label{11}
\vec{f}(\rho,\T) = \vec{0}.
\end{equation}
The above is a generalization of the explicit Euler constitutive relation to the class of implicit relations. The class of models defined through \eqref{9} form a sub-class of the implicit models considered by Rajagopal (2003). The representation \eqref{7} simplifies to
\begin{equation}
    \label{12}
    \alpha_1 \I + \alpha_2 \T + \alpha_4 \T^2 = \vec{0},
\end{equation}
where the material moduli $\alpha_i$, $i=1,2,4$, depend on the density $\rho$ and the invariants $\tr\T$, $\tr\T^2$ and $\tr\T^3$. 
Notice the difference between the representation \eqref{12} that the stress satisfies and the expression for the stress in a Euler fluid. Even if $\alpha_4=0$, the representation \eqref{12} is different and is more general than the constitutive relation for a Euler fluid as the material moduli $\alpha_i$, $i=1,2$ depend on the density  and the invariants $\tr\T$, $\tr\T^2$ and $\tr\T^3$. That is, even in the case of $\alpha_4 = 0$, the equation \eqref{12} reduces to
\begin{equation}
    \label{13}
    \alpha_1 (\rho, \tr\T, \tr\T^2, \tr\T^3)\I + \alpha_2 (\rho, \tr\T, \tr\T^2, \tr\T^3) \T  = \vec{0},
\end{equation}
that one may not be able to solve to obtain an explicit expression for the stress in term of the density. We note that the Euler fluid is a special case of \eqref{13} and corresponds to $\alpha_1$ being a function of the density $\rho$ and $\alpha_2$ is a constant. The question is whether the additional structure in \eqref{13} allows us to describe additional natural phenomena or experiments than the Euler fluid. Before we proceed to do this, we shall discuss the third purpose of the paper, the challenge offered by these new implicit constitutive relations in mathematical and numerical analysis.

We now record the governing equations that we need to consider. The balance of mass is given by
\begin{equation}
    \label{14} 
    \frac{\partial \rho}{\partial t} + \diver(\rho\vec{v}) = 0,
\end{equation}
and the balance of linear momentum is given by
\begin{equation}
    \label{15} 
    \rho\frac{\textrm{d} \vec{v}}{\textrm{dt}} = \diver\T + \rho\vec{b},
\end{equation}
where $\vec{b}$ is the specific body force. We also assume that the Cauchy stress $\T$ is symmetric and the balance of angular momentum is fulfilled. Thus, we now need to solve the system of equations \eqref{12}, \eqref{14}, and \eqref{15} (ten scalar equations) simultaneously for the unknowns $\rho$, $\vec{v}$, $\T$ (ten scalar unknowns). Unlike the situation in the case of a classical Euler fluid, we cannot substitute the expression for the stress $\T$  into the balance of linear momentum and get an equation that relates the gradient of the pressure and the velocity. Here, we need to solve the balance laws in conjunction with the constitutive equation \eqref{10}.      
      
In the case of the generalization of the Korteweg model, we need to solve equations \eqref{7}, \eqref{14}, and \eqref{15} simultaneously for the unknowns $\rho$, $\vec{v}$, $\T$. Once again, ten partial differential equations for ten unknowns. Such systems can present very interesting challenges to study issues such as existence, uniqueness and stability.

\section{A simple boundary value problem that can be described by an infinite sub-class of constitutive relations belonging to the class \eqref{12}}\label{Sec3}

Let us consider a very simple problem within the context of the constitutive relation \eqref{12}. Let us consider a half-space $\{(x,y,z); x,z\in (-\infty,\infty) \textrm{ and } y\in (-\infty, 0)\}$ filled with a fluid described by \eqref{12} at rest under the action of gravity whose density does not change with time. Let us assume that the free surface is defined by $y=0$. Since the fluid is static, the velocity $\vec{v} =\vec{0}$. Then, since the density does not change with time, \eqref{14} is automatically satisfied. Thus, we need to find a stress $\T$ that satisfies \eqref{12} and \eqref{15} simultaneously\footnote{As $\vec{v}=\vec{0}$, $\frac{\textrm{d}\vec{v}}{\textrm{dt}} = \vec{0}$ in \eqref{15}.}. We shall now investigate whether a specific flow field can be satisfied by a sub-class of constitutive relations \eqref{12} and meet the requirement of \eqref{15}.

Let us suppose the body force (gravity) is given by
\begin{equation}
    \label{16}
    \vec{b} = -g\vec{j} = (0, -g, 0), 
\end{equation}
where $g$ is the accelaration due  to gravity and $\vec{j}$ denotes the unit vector in the $y$-coordinate direction. Furthermore, we shall appeal to the semi-inverse method to seek a solution for the stress field of the form
\begin{equation}
    \label{17}
    \T = -\varphi(y) \I.  
\end{equation}
It then immediately follows from \eqref{15} that
\begin{equation}
    \label{18}
    \varphi(y) =  -  \int g \rho(s) \d s \qquad \textrm{ or } \quad \varphi(y) = \varphi(0) + \int_{y}^0 g\rho(s) \d s.   
\end{equation} 
We now must verify if this solution satisfies \eqref{12}. Since 
\begin{equation}
    \label{19}
    \tr\T = -3\varphi, \quad \tr \T^2 = 3 \varphi^2 \quad \textrm{ and } \tr \T^3 = -3 \varphi^3,   
\end{equation}
it follows from \eqref{12} that we need to satisfy 
\begin{equation}
    \label{20}
    \hat{\alpha}_1 (\rho, \varphi) - \hat{\alpha}_2 (\rho, \varphi) \varphi + \hat{\alpha}_4 (\rho, \varphi) \varphi^2 = 0\,
\end{equation}
where $\hat\alpha_i(\rho, \varphi)$ denotes $\alpha_i(\rho, 3\varphi, - 3\varphi^2, 3\varphi^3)$, $i=1,2,4$. Thus, a solution is not possible for all constitutive relations of the form \eqref{12}, solutions are only possible if the constitutive relations are such that equation \eqref{20} is satisfied when $\varphi$ is given by \eqref{18}.

The Euler fluid is a special case of the implicit constitutive relation. Let us consider the same problem but within the context of a Euler fluid whose constitutive relation takes the form
\begin{equation}
    \label{21}
    \T = -p(\rho)\I, 
\end{equation}       
and furthermore, suppose that $p$ takes the special form
\begin{equation}\label{22}
    p(\rho)=C\rho,
\end{equation} 
where $C$ is a constant, that is we are considering an ideal gas.     
                                                                      
Let us consider such a fluid is at rest under the action of gravity, and let us suppose that the density is just a function of $y$. It then follows from the balance of linear momentum that
\begin{equation}\label{23}
    \rho = K\exp\left( - \frac{g}{C} y\right),
\end{equation} 
where the constant $K$ is given by $K=\rho(0)$. We however do not know the value of the density at $y=0$. It follows from (22) and (23) that
\begin{equation}\label{24}
    p = KC\exp\left( - \frac{g}{C} y\right).
\end{equation} 
In order to determine the constant $K$, we could possibly evaluate the pressure at $y=0$, and this is what is usually done. Different choices for $p$ as a function of $\rho$ would lead to different solutions for the density and the pressure. 

The interesting question is whether one can find the solution of the form \eqref{23} for a constitutive relation belonging to \eqref{12} other than the Euler constitutive relation, that is can one find a $\varphi$ such that the stress is given by \eqref{17} such that \eqref{12} and \eqref{23} are met. After inserting \eqref{23} into \eqref{20}, one obtains
\begin{equation}
    \label{25}
    \hat{\alpha}_1 \left(K\exp\left( - \frac{g}{C} y\right), \varphi\right) - \hat{\alpha}_2 \left(K\exp\left( - \frac{g}{C} y\right), \varphi\right) \varphi + \hat{\alpha}_4 \left(K\exp\left( - \frac{g}{C} y\right), \varphi\right) \varphi^2 = 0,
\end{equation}
and the question reduces to whether there exist $\alpha_i$, $i=1,2,3$ such that we can find a solution to \eqref{25} other than the classical solution for the Euler fluid. As an example, setting 
\begin{equation}
    \label{26}
    \begin{split}
    \hat\alpha_1\left(K\exp\left( - \frac{g}{C} y\right), \varphi\right) &= A\varphi \exp\left( - \frac{g}{C} y\right), \\ \hat\alpha_2\left(K\exp\left( - \frac{g}{C} y\right), \varphi\right) &= A\exp\left( - \frac{g}{C} y\right), \\ 
    \hat\alpha_4\left(K\exp\left( - \frac{g}{C} y\right), \varphi\right) &=0,
    \end{split}
\end{equation}
we fulfil \eqref{25}. So would  
\begin{equation}
    \label{27}
    \begin{split}
    \hat\alpha_1\left(K\exp\left( - \frac{g}{C} y\right), \varphi\right) &= A\varphi^2 \exp\left( - \frac{g}{C} y\right), \\ \hat\alpha_2\left(K\exp\left( - \frac{g}{C} y\right), \varphi\right) &= A\varphi\exp\left( - \frac{g}{C} y\right), \\ 
    \hat\alpha_4\left(K\exp\left( - \frac{g}{C} y\right), \varphi\right) &=0.
    \end{split}
\end{equation}
In fact, infinity of constitutive relations would satisfy the requirement. The fact that more than one model belonging to the class of implicit functions can describe a particular phenomenon has important consequences. If one were to restrict oneself to the class of classical Euler models, then the solution \eqref{23} will be possible in only the Euler fluid defined through \eqref{22}. On the other hand, infinity of fluids described by the class of implicit models can support the same solution. Thus, while two models belonging to the class of implicit relations may both be capable of describing a particular phenomenon, each of them might be capable of describing different observations (phenomenon/experimental result). However, having described a particular result a specific Euler fluid may be incapable of describing some other observation. Thus, by considering a sequence of boundary value problems, we can determine the constitutive relations that belong to \eqref{12} that can best explain all of them.

To recapitulate, as the result cannot be overemphasized, the class of implicit constitutive relations \eqref{12} allows one to have many constitutive relations that describe the result \eqref{23}. Given a set of observations, we can systematically cull the class of constitutive relations to arrive at a sub-class which best describes the set of observations. This culling process may lead to one or for that matter no constitutive relation that can describe the set of observations. A good example of the latter is the class of observation of turbulent flows. To date, no adequate theory has been found that can describe well all the observed turbulent phenomena.

\bibliographystyle{amsplain}
\bibliography{biblio}

\end{document}